\documentclass[12pt]{amsart}%
\usepackage{amsfonts}
\usepackage{amsmath}
\usepackage{amssymb}
\usepackage{graphicx}%
\makeatletter
\@addtoreset{equation}{section}
\makeatother

\newtheorem{theorem}{Theorem}[section]

\newtheorem{corollary}[theorem]{Corollary}

\newtheorem{lemma}[theorem]{Lemma}

\newtheorem{remark}[theorem]{Remark}

\makeatletter
\def\@makefnmark{}
\makeatother

\begin{document}

\title[Stein-Weiss inequalities with the fractional Poisson kernel]{Stein-Weiss inequalities with the fractional Poisson kernel}

\author{Lu Chen}

\address{School of Mathematical and statistics\\
Beijing Institute of Technology, Beijing 100081, P. R. China P. R. China}
\email{chenlu5818804@163.com}

\author{Zhao Liu}
\address{School of Mathematics and Computer Science, Jiangxi Science Technology Normal University, Nanchang 330038, People's Republic of China}
\email{liuzhao2008tj@sina.com}

\author{Guozhen Lu }
\address{Department of Mathematics\\
University of Connecticut\\
Storrs, CT 06269, USA}
\email{guozhen.lu@uconn.edu}

\author{Chunxia Tao}
\address{School  of Mathematical Sciences\\
 Beijing Normal  University\\
  Beijing 100875, China}
\email{taochunxia@mail.bnu.edu.cn}

\thanks{The first and fourth authors were partly supported by a grant from the NNSF of China (No.11371056), the second author was partly
supported by a grant from the NNSF of China (No.11801237), the third author was partly supported by a grant from the Simons Foundation.
}

\begin{abstract}

In this paper, we establish the following Stein-Weiss inequality with the fractional Poisson kernel (see Theorem 1.1):
\begin{equation}\label{int1}
\int_{\mathbb{R}^n_{+}}\int_{\partial\mathbb{R}^n_{+}}|\xi|^{-\alpha}f(\xi)P(x,\xi,\gamma)g(x)|x|^{-\beta}d\xi dx\leq C_{n,\alpha,\beta,p,q'}\|g\|_{L^{q'}(\mathbb{R}^n_{+})}\|f\|_{L^p(\partial \mathbb{R}^{n}_{+})},
\end{equation}
where $P(x,\xi,\gamma)=\frac{x_n}{(|x'-\xi|^2+x_n^2)^{\frac{n+2-\gamma}{2}}}$, $2\le \gamma<n$, $f\in L^{p}(\partial\mathbb{R}^n_{+})$, $g\in L^{q'}(\mathbb{R}^n_{+})$ and $p,\ q'\in (1,\infty)$ and satisfy $\frac{n-1}{n}\frac{1}{p}+\frac{1}{q'}+\frac{\alpha+\beta+2-\gamma}{n}=1$. Then we prove that there exist extremals for the Stein-Weiss inequality \eqref{int1}   and  the extremals  must be radially decreasing about the origin (see Theorem 1.5). We also
 provide the regularity and asymptotic estimates of positive solutions to the integral
systems which are the Euler-Lagrange equations of the extremals to the Stein-Weiss inequality \eqref{int1} with the fractional Poisson kernel (see Theorems 1.7 and 1.8).
Our result is inspired by the work of Hang, Wang and Yan \cite{HWY} where the Hardy-Littlewood-Sobolev type inequality was first established
when $\gamma=2$ and $\alpha=\beta=0$ (see \eqref{PI}). The proof of the Stein-Weiss inequality (\ref{int1}) with the fractional Poisson kernel in this paper uses our recent work on the Hardy-Littlewood-Sobolev inequality with the fractional Poisson kernel \cite{CLT} and the present paper is a further study in this direction.

\end{abstract}
\maketitle {\small {\bf Keywords:} Existence of extremal functions; Stein-Weiss inequality; Poisson kernel; Hardy inequality in high dimensions. \\

}

\section{Introduction}
The study and understanding of various kinds of weighted integral inequalities has attracted a great attention of many people due to the importance of such inequalities in applications to problems in harmonic analysis and partial differential equations. Generally, these inequalities play a key role in establishing the existence and radial symmetry results for certain non-linear equations. Let us first recall the Stein-Weiss inequality.
\medskip

The well-known Stein-Weiss inequality which was established by Stein and Weiss in \cite{Stein} states that
\begin{equation}\label{SW}
\int_{\mathbb{R}^n}\int_{\mathbb{R}^n}|x|^{-\alpha}|x-y|^{-\lambda} f(x)g(y)|y|^{-\beta} dxdy\leq C_{n,\alpha,\beta,p,q'}\|f\|_{L^{q'}(\mathbb{R}^n)}\|g\|_{L^p (\mathbb{R}^n)},
\end{equation}
where $1<p$, $q'<\infty$, $\alpha$, $\beta$ and $\lambda$ satisfy the following conditions,
$$\frac{1}{q'}+\frac{1}{p}+\frac{\alpha+\beta+\lambda}{n}=2,\ \ \frac{1}{q'}+\frac{1}{p}\geq 1,$$
$$\alpha+\beta\geq0,\ \ \alpha<\frac{n}{q},\ \ \beta<\frac{n}{p'},\ \ 0<\lambda<n.$$
(see also an alternative proof of establishing the Stein-Weiss inequalities  recently found in \cite{HLZ}  by using conditions on weights to guarantee the weighted boundedness of fractional integrals given in \cite{SawyerWheeden} and such a method also applies to establish the Stein-Weiss inequalities on the Heisenberg groups.)
Lieb \cite{Lieb} used the method based on symmetrization argument and the Riesz rearrangement     to establish the existence of extremals for the inequality \eqref{SW} in the case $p<q$ and $\alpha$, $\beta \geq 0$. Furthermore, in the case of $p=q$,   the extremals can't be expected to exist (see Lieb \cite{Lieb} and also  Herbst \cite{Herbst} for the case $\lambda=n-1$, $p=q=2$, $\alpha=0$, $\beta=1$). In the case of $p=q$, Beckner \cite{B2,B3} obtained the sharp constant of the Stein-Weiss inequalities \eqref{SW} by establishing an equivalent formulation as a convolution estimate on the product manifold $\mathbb{R}^{+}\times \mathbb{S}^{n-1}$. The precise estimate of the sharp constant of the Stein-Weiss inequalities for the case of $p\neq q$ was also established in \cite{B3}. For more results about proving precise estimates for Stein-Weiss functionals in conjunction with the study of Housdorff-Young and Pitt's type inequalities and their multilinear versions, we refer the reader to the works of Beckner  \cite{B0,B1, B4, B5, B6, B7}.
We note that the existence of extremal functions for the Stein-Weiss inequalities in the case $p<q$ under the assumption $\alpha+\beta \geq 0$ has been established by Chen, Lu and Tao \cite{CLT1}, which extends Lieb's result under the stronger assumption that $\alpha\geq 0$ and $\beta \geq 0$, using the concentration-compactness of Lions (\cite{Lions1, Lions2}).

\medskip

In the special case of $\alpha=\beta=0$, the Stein-Weiss inequality \eqref{SW} becomes the following Hardy-Littlewood-Sobolev  (HLS) inequality (see \cite{Hardy, Sobolev}),
\begin{equation}\label{HL1}
\int_{\mathbb{R}^n}\int_{\mathbb{R}^n}|x-y|^{-\lambda} f(x)g(y)dxdy\leq C_{n,p,q'}\|f\|_{L^{q'}(\mathbb{R}^n)}\|g\|_{L^p(\mathbb{R}^n)},
\end{equation}
where $1<q', p<\infty, 0<\lambda<n$ and $\frac{1}{q'}+\frac{1}{p}+\frac{\lambda}{n}=2$.
When one of $p$ and $q'$ equals 2 or $p = q'$, Lieb \cite{Lieb} obtained the sharp constants of this inequality.  
By using the competing symmetry method, Carlen and Loss \cite{CarlenLoss}  provided a different proof from Lieb's of the sharp constants and extremal functions in the diagonal case $p=q^{'}=\frac{2n}{2n-\lambda}$ and Frank and Lieb \cite{FL1} offered a new proof using the reflection positivity of inversions in spheres in the special diagonal case. Carlen, Carillo and Loss gave a simple proof of the  sharp Hardy-Littlewood-Sobolev
inequality when $\lambda=n-2$ for $n\geq 3$ via a monotone
flow governed by the fast diffusion equation \cite{CCL}.
Frank and Lieb \cite{FL2} further employed a rearrangement-free technique developed in \cite{Lieb1} to recapture  the best constant of inequality \eqref{HL1}. When $q'=p =\frac{2n}{2n-\lambda}$, Euler-Lagrange equation of the extremals to the HLS inequality is a conformal invariant integral equation. Using the method of moving plane or moving-sphere in integral forms (see \cite{CLO, Li}), one can classify the positive solutions to this integral equation. The HLS and Stein-Weiss inequalities also have many applications in partial differential equations. One can also see \cite{B2, CJLL, CL, CLL,LZ} and the references therein.
\vskip0.3cm

Dou \cite{Dou0} established a double weighted HLS inequality on upper half space $\mathbb{R}^{n}_{+}$ :
\begin{equation}\label{SWU}
\int_{\mathbb{R}^n_{+}}\int_{\partial\mathbb{R}^n_{+}}|x|^{-\alpha}|x-y|^{-\lambda} f(x)g(y)|y|^{-\beta} dxdy\leq C_{n,\alpha,\beta,p,q'}\|f\|_{L^p(\partial \mathbb{R}^n_{+})}\|g\|_{L^{q'}(\mathbb{R}^{n}_{+})},
\end{equation}
where  $p$, $q'$, $\alpha$, $\beta$ and $\lambda$ satisfy the following conditions:
$$\frac{1}{q'}+\frac{n-1}{np} +\frac{\alpha+\beta+\lambda+1}{n}=2,\ \ \frac{1}{q'}+\frac{1}{p}\geq 1,$$
$$\alpha+\beta\geq0,\ \ \alpha<\frac{n-1}{p'},\ \ \beta<\frac{n}{q},\ \ 0<\lambda<n-1.$$
In the case $\alpha$, $\beta\geq 0$ and $p<q$, he also discussed the existence of extremals for this inequality through the rearrangement inequality.
If we set $\alpha=\beta=0$ in the inequality \eqref{SWU}, then this inequality reduces to the HLS inequality on upper half space which was studied by Dou and Zhu \cite{Dou1}.

\medskip

Recently, Chen, Lu and Tao \cite{CLT} considered the sharp HLS inequality with the fractional Poisson kernel which can be seen as a extension of inequality \eqref{HL1},
\begin{equation}\label{PHLS1}
\int_{\mathbb{R}^n_{+}}\int_{\partial\mathbb{R}^n_{+}}f(\xi)P(x,\xi,\gamma)g(x)d\xi dx\leq C_{n,\alpha,p,q'}\|g\|_{L^{q'}(\mathbb{R}^n_{+})}\|f\|_{L^p(\partial \mathbb{R}^{n}_{+})},
\end{equation}
where $1<p, q'<\infty$ and $2\leq \gamma<n$, satisfying $\frac{n-1}{n}\frac{1}{p}+\frac{1}{q'}+\frac{2-\gamma}{n}=1$ and $P(x,\xi,\gamma)$ is the fractional Poisson kernel which will be defined below.
They employed the rearrangement inequality and Lorentz interpolation to obtain the existence of extremals to the inequality \eqref{PHLS1}. Furthermore,
the radial symmetry of extremals was also established through the method of moving-plane in integral forms. We also note that
$P(x,\xi,\gamma)$ can be seen as the fundamental solution of the fractional Laplacian operator. In fact, when $\gamma=2$, $P(x,\xi,\gamma)$
is actually the classical Poisson kernel corresponding to the fundamental solution of the Laplacian operator. Inequality \eqref{PHLS1} for $\gamma=2$ is
equivalent to the following integral inequality with the Poisson kernel which was first established by Hang, Yang and Wang in \cite{HWY},
\begin{equation}\label{PI}
\|\int_{\partial\mathbb{R}^n_{+}}P(x,\xi,2)f(\xi)d\xi \|_{L^{q}(\mathbb{R}^n_{+})}\leq C(n,p)\|f\|_{L^{p}(\partial\mathbb{R}^n_{+})},
\end{equation}
where $P(x,\xi,2)=c_n\frac{x_n}{(|x'-\xi|^2+x_n^2)^{\frac{n}{2}}}$ and $q=\frac{np}{n-1}$.
\medskip

In the spirit of the aforementioned  works,  we are interested to investigate whether the HLS type inequality \eqref{PHLS1} can be extended to a weighted HLS inequality with the fractional Poisson kernel. Furthermore, we would like to know if  such an inequality have an extremal function for all the indices.

In  the present paper, we provide positive answers to these questions. We first establish Stein-Weiss inequality with the fractional Poisson kernel.
\medskip
\begin{theorem}\label{theorem1}
For $n\geq 3$, $1<p< \infty$, $1< q'<\infty$, $\frac{n-1}{n}\frac{1}{p}+\frac{1}{q'}\geq1$, $2\leq \gamma<n$, $\alpha<\frac{n-1}{p'}$,  $\beta<\frac{n}{q}+1$, $\alpha+\beta\geq0$, satisfying
\begin{equation}\label{tiaojian}
\frac{n-1}{n}\frac{1}{p}+\frac{1}{q'}+\frac{\alpha+\beta+2-\gamma}{n}=1,
\end{equation}
then there exists some constant $C_{n,\alpha,p,q'}>0$ such that for any functions  $f\in L^p(\partial\mathbb{R}^n_{+})$ and $g\in L^{q'}(\mathbb{R}^n_{+})$, there holds
\begin{equation}\label{PHLS}
\int_{\mathbb{R}^n_{+}}\int_{\partial\mathbb{R}^n_{+}}|\xi|^{-\alpha}f(\xi)P(x,\xi,\gamma)g(x)|x|^{-\beta}d\xi dx\leq C_{n,\alpha,\beta,p,q'}\|g\|_{L^{q'}(\mathbb{R}^n_{+})}\|f\|_{L^p(\partial \mathbb{R}^{n}_{+})},
\end{equation}
where $P(x,\xi,\gamma)=\frac{x_n}{(|x'-\xi|^2+x_n^2)^{\frac{n+2-\gamma}{2}}}$ with $x=(x',x_n)\in \mathbb{R}^{n-1}\times \mathbb{R}^{+}$.
\end{theorem}
\begin{remark}
When $\alpha=\beta=0$, $p=\frac{2(n-1)}{n+\alpha-2}$ and $q'=\frac{2n}{n+\alpha}$, Dou, Guo and Zhu \cite{DGZ} applied the methods based on
conformal transformation and moving sphere in integral forms to obtain extremal functions of the
inequality \eqref{PHLS} and computed the sharp constant.
\end{remark}

Define the double weighted operator
\begin{equation}\label{fractionaloperator}
V(f)(x)=\int_{\mathbb{\partial R}^n_{+}} |\xi|^{-\alpha}P(x,\xi,\gamma) f(\xi)|x|^{-\beta}d\xi,\ W(g)(\xi)=\int_{\mathbb{R}^n_{+}} |\xi|^{-\alpha}P(x,\xi,\gamma) g(x)|x|^{-\beta}dx.
\end{equation}
By duality, we can see that the inequality \eqref{PHLS} is equivalent to the following weighted inequality
\begin{equation}\label{we1}
\|V(f)\|_{L^q(\mathbb{R}^n_{+})}\leq C_{n,\alpha,\beta,p,q'}\|f\|_{L^p(\partial \mathbb{R}^{n}_{+})},\ \|W(g)\|_{L^{p'}(\partial\mathbb{R}^n_{+})}
\leq  C_{n,\alpha,\beta,p,q'} \|g\|_{L^{q'}(\mathbb{R}^n_{+})}.
\end{equation}
where $\frac{1}{q}=\frac{n-1}{n}\frac{1}{p}+\frac{\alpha+\beta+2-\gamma}{n}$ and $\frac{1}{p'}=\frac{n}{n-1}\frac{1}{q'}+\frac{\alpha+\beta+1-\gamma}{n-1}$.

By the Marcinkiewicz interpolation theorem (see \cite{Stein2}), we immediately obtain the following stronger results involving the Lorentz norm $L^{p, q}(\partial \mathbb{R}^n_+)$.
\begin{corollary}\label{corollary2}
Assume that  $p$, $q$, $\alpha$, $\beta$ and $\gamma$ satisfy the conditions of Theorem \ref{theorem1}, then there holds
\begin{equation}\label{PHLS4}
\|V(f)\|_{L^{q}(\mathbb{R}^n_{+})}\leq C_{n,\alpha,\beta,p,q'}\|f\|_{L^{p,q}(\partial \mathbb{R}^{n}_{+})}.
\end{equation}
\end{corollary}

When $\gamma=2$, we can conclude the following Stein-Weiss inequality with the Poisson kernel from Theorem \ref{theorem1}.
\begin{corollary}
For $n\geq 3$, $1<p< \infty$, $1< q'<\infty$, $\alpha<\frac{n-1}{p'}$,  $\beta<\frac{n}{q}+1$, $\alpha+\beta=0$, satisfying
\begin{equation}\label{tiaojian}
\frac{n-1}{n}\frac{1}{p}+\frac{1}{q'}=1,
\end{equation}
then there exists some constant $C_{n,\alpha,p,q'}>0$ such that for any functions $f\in L^p(\partial\mathbb{R}^n_{+})$ and $g\in L^{q'}(\mathbb{R}^n_{+})$, there holds
\begin{equation*}\label{PHLS2}
\int_{\mathbb{R}^n_{+}}\int_{\partial\mathbb{R}^n_{+}}|\xi|^{-\alpha}f(\xi)P(x,\xi,2)g(x)|x|^{-\beta}d\xi dx\leq C_{n,\alpha,\beta,p,q'}\|g\|_{L^{q'}(\mathbb{R}^n_{+})}\|f\|_{L^p(\partial \mathbb{R}^{n}_{+})}.
\end{equation*}
\end{corollary}

Next, we start to establish the existence of a maximizing pair of functions $(f(\xi)$, $g(x))$ giving equality in \eqref{PHLS}, which means
finding $f(\xi) \in L^{p}(\partial \mathbb{R}^{n}_{+})$ with $\|f\|_{L^p(\partial \mathbb{R}^{n}_{+})}=1$ and  $\frac{1}{q}=\frac{n-1}{n}\frac{1}{p}+\frac{\alpha+\beta+2-\gamma}{n}$
such that
$$\|V(f)\|_{L^q(\mathbb{R}^n_{+})}=\sup\{\|V(h)\|_{L^q(\mathbb{R}^n_{+})} : h\geq 0 ,\  \|h\|_{L^p(\partial \mathbb{R}^{n}_{+})}=1\}=C_{n,\alpha,\beta,p,q'}.$$
In the following result, we prove such an extremal function for \eqref{PHLS} indeed exists.
\begin{theorem}\label{theorem2}
For $n\geq 3$, $1<p< \infty$, $1< q'<\infty$, $\frac{n-1}{n}\frac{1}{p}+\frac{1}{q'}\geq1$, $2\leq \gamma<n$, $\alpha<\frac{n-1}{p'}$,  $\beta<\frac{n}{q}+1$, $\alpha$, $\beta\geq 0$ satisfying $$\frac{1}{q}=\frac{n-1}{n}\frac{1}{p}+\frac{\alpha+\beta+2-\gamma}{n},$$
there exists some nonnegative function $f\in L^p(\mathbb{\partial R}^n_{+})$ satisfying $\|f\|_{L^p(\partial\mathbb{R}^n_{+})}=1$ and $\|V(f)\|_{L^q(\mathbb{R}^n_{+})}=C_{n,\alpha,\beta,p,q'}$. Furthermore, if $(f(\xi), g(x))$ is a pair of maximizer of the inequality \eqref{PHLS}, then
$f(\xi)$ is radially symmetric and monotone decreasing about the origin and $g(x)=c_0V(f)(x)$ for some $c_0>0$.
\end{theorem}

\begin{remark}
Given $h$ a measurable function on $\partial \mathbb{R}^n_+$, let $h^*$ be the decreasing rearrangement of $h$.
According to the Brascamp-Lieb-Luttinger rearrangement inequality (see \cite{BLL}), we have
$$\|h^{*}\|_{L^{p}(\partial\mathbb{R}^n_{+})}=\|h\|_{L^{p}(\partial\mathbb{R}^n_{+})},\ \
\|V(h)\|_{L^{q}(\mathbb{R}^n_{+})}\leq\|V(h^{*})\|_{L^{q}(\mathbb{R}^n_{+})}.$$
Hence if $(f, g)$ is a pair of maximizer of the inequality \eqref{PHLS}, then we must have
$$\|V(f)\|_{L^{q}(\mathbb{R}^n_{+})}=\|V(f^{*})\|_{L^{q}(\mathbb{R}^n_{+})},$$
which implies that $f$ is radially symmetric and monotone decreasing about the origin.
Noting that $(f, g)$ is a pair of maximizer of the inequality \eqref{PHLS}, using the H\"{o}lder inequality, we derive that
\begin{equation}\begin{split}
C_{n,\alpha,\beta,p,q'}\|g\|_{L^{q'}(\mathbb{R}^n_{+})}\|f\|_{L^p(\partial \mathbb{R}^{n}_{+})}&=\int_{\mathbb{R}^n_{+}}V(f)(x)g(x)dx\\
&\leq \|V(f)\|_{L^{q}(\mathbb{R}^n_{+})}\|g\|_{L^{q'}(\mathbb{R}^n_{+})}\\
&\leq C_{n,\alpha,\beta,p,q'}\|g\|_{L^{q'}(\mathbb{R}^n_{+})}\|f\|_{L^p(\partial \mathbb{R}^{n}_{+})}.
\end{split}\end{equation}
Hence the H\"{o}lder inequality $$\int_{\mathbb{R}^n_{+}}V(f)(x)g(x)dx\leq \|V(f)\|_{L^{q}(\mathbb{R}^n_{+})}\|g\|_{L^{q'}(\mathbb{R}^n_{+})}$$
must be an equality, which implies that $g(x)=c_0V(f)(x)$ for some $c_0>0$.
\end{remark}
Since we have established the existence of extremals to the inequalities \eqref{PHLS},
we are naturally concerned with some properties such as the regularity and asymptotic behavior of the extremals.
To this purpose, by maximizing the functional
\begin{equation}\label{functional}
J(f,g)=\int_{\mathbb{R}^n_{+}}\int_{\partial\mathbb{R}^n_{+}}|\xi|^{-\alpha}f(\xi)|P(x,\xi,\gamma)g(x)|x|^{-\beta}d\xi dx
\end{equation}
under the constraint $\|f\|_{L^{p}(\partial\mathbb{R}^n_{+})}=\|g\|_{L^{q'}(\mathbb{R}^n_{+})}=1$ and Euler-Lagrange multiplier theorem, we can derive the following double weighted integral system associated with the fractional Poisson kernel.
\begin{equation}\label{eu}\begin{cases}
J(f,g)f(\xi)^{p-1}=\int_{\mathbb{R}^n_{+}}|\xi|^{-\alpha}P(x,\xi,\gamma)g(x)|x|^{-\beta}dx,\ \ \xi \in \partial\mathbb{R}^{n}_{+},\\
J(f,g)g(x)^{q'-1}=\int_{\partial\mathbb{R}^n_{+}}|x|^{-\beta}P(x,\xi,\gamma)f(\xi)|\xi|^{-\alpha} d\xi, \ \ x\in \mathbb{R}^{n}_{+}.
\end{cases}\end{equation}
\medskip

Let $u=c_1f^{p-1}, v=c_2g^{q'-1}, \frac{1}{p-1}=p_0$ and $\frac{1}{q'-1}=q_0$ and pick two suitable
constants $c_1$ and $c_2$, then system \eqref{eu} is simplified as
\begin{equation}\label{system3}\begin{cases}
u(\xi)=\int_{\mathbb{R}^n_{+}} |\xi|^{-\alpha}P(x,\xi,\gamma)v^{q_0}(x)|x|^{-\beta}dx,\ \  \xi\in \partial\mathbb{R}^n_{+},\\
v(x)=\int_{\partial\mathbb{R}^n_{+}}|x|^{-\beta}P(x,\xi,\gamma)u^{p_0}(\xi)|\xi|^{-\alpha} d\xi,\ \ x\in \mathbb{R}^n_{+},
\end{cases}\end{equation}
where $\alpha$, $\beta\geq 0$, $p_0$ and $q_0$ satisfy $\frac{n-1}{n}\frac{1}{p_0+1}+\frac{1}{q_0+1}=\frac{n+\alpha+\beta+1-\gamma}{n}.$
\medskip

Next, we give the regularity estimates and asymptotic behaviors of the solutions to the integral system \eqref{system3}.

\begin{theorem}\label{thmx1}
Assume that $(u,v)\in L^{p_0+1}(\partial\mathbb{R}^n_{+})\times L^{q_0+1}(\mathbb{R}^n_{+})$ is a pair of positive solutions of the integral system \eqref{system3}, then $(u,v)\in L^{r}(\partial\mathbb{R}^n_{+})\times L^{s}(\mathbb{R}^n_{+})$ for all $r$ and $s$ such that
$$\frac{1}{r}\in \Big(\frac{\alpha}{n-1}, \frac{n+2-\gamma+\alpha}{n-1}\Big)\cap \Big(\frac{1}{p_0+1}-\frac{n}{n-1}\frac{1}{q_0+1}+\frac{\beta-1}{n-1}, \frac{1}{p_0+1}-\frac{n}{n-1}\frac{1}{q_0+1}+\frac{n+1-\gamma+\beta}{n-1}\Big)$$
and
$$\frac{1}{s}\in \Big(\frac{\beta-1}{n}, \frac{n+1-\gamma+\beta}{n}\Big)\cap \Big(\frac{1}{q_0+1}-\frac{n-1}{n}\frac{1}{p_0+1}+\frac{\alpha}{n}, \frac{1}{q_0+1}-\frac{n-1}{n}\frac{1}{p_0+1}+\frac{n+2-\gamma+\alpha}{n}\Big).$$
\end{theorem}

\begin{theorem}\label{thmx2}
Assume that $(u,v)\in L^{p_0+1}(\partial\mathbb{R}^n_{+})\times L^{q_0+1}(\mathbb{R}^n_{+})$ is a pair of positive solutions of the integral system \eqref{system3}, suppose that $p_0$, $q_0>1$, then there holds $$\lim\limits_{|\xi|\rightarrow 0}u(\xi)|\xi|^{\alpha}=
\int_{\mathbb{R}^n_{+}}\frac{v^{q_0}(x)x_n}{|x|^{n+2-\gamma+\beta}}dx$$ if
$\frac{1}{q_0}-\frac{n+1+\beta-\gamma}{q_0n}>\frac{\beta-1}{n}$. Similarly,
$$\lim\limits_{|x|\rightarrow 0}\frac{v(x)|x|^{\beta}}{x_n}=\int_{\partial\mathbb{R}^n_{+}}\frac{|u|^{p_0}(\xi)}{|\xi|^{n+2-\gamma+\alpha}}d\xi$$
if $\frac{1}{p_0}-\frac{n+2+\alpha-\gamma}{p_0(n-1)}>\frac{\alpha}{n-1}$.
\end{theorem}

Finally, we are also interested in studying the following single weighted integral system with the fractional Poisson kernel
\begin{equation}\label{int}\begin{cases}
u(\xi)=\int_{\mathbb{R}^n_{+}} P(x,\xi,\gamma) |x|^{-\beta} v^{q_0}(x)dx,\ \ \xi\in \partial \mathbb{R}^n_{+},\\
v(x)=\int_{\partial\mathbb{R}^n} P(x,\xi,\gamma) |\xi|^{-\alpha} u^{p_0}(\xi)d\xi, \ \ x\in \mathbb{R}^n_{+}.
\end{cases}\end{equation}
With the help of the Pohozaev identity, we establish the following necessary condition.
\begin{theorem}\label{theorem4}
For $2<\gamma<n$, $0< p_0<\infty$, $0< q_0<\infty$, suppose that there exists a pair of $C^1$ positive solutions $(u,v)\in L^{p_0+1}(|\xi|^{-\alpha}d\xi, \partial\mathbb{R}^n_{+}) \times L^{q_0+1}(|x|^{-\beta}dx, \mathbb{R}^n_{+})$ satisfying the integral system \eqref{int}, then the following balance condition must hold:
$$\frac{n-1-\alpha}{p_0+1}+\frac{n-\beta}{q_0+1}=n+1-\gamma.$$
\end{theorem}
\medskip

This paper is organized as follows. In Section 2, we employ the weighted Hardy inequality in high dimensions and the HLS inequality with the fractional Poisson kernel inequality to establish the Stein-Weiss inequality with the fractional Poisson kernel \eqref{PHLS}. In Section 3, we discuss the existence of extremals to the inequality \eqref{PHLS} through the rearrangement argument. Sections 4 and 5 are devoted to the regularity and asymptotic behavior to the double weighted integral system associated with the fractional Poisson kernel.

\medskip

{\bf Acknowledgement:} The authors wish to thank the referee for his many helpful comments on the paper which have improved the exposition and for many
constructive suggestions to improve the paper. In the earlier version of the paper, the radial symmetry of the extremal functions was established by using the moving plane method in integral forms. It is pointed out by the referee that the radial symmetry of the extremal functions
of the Stein-Weiss inequality follows quickly from the Brascamp-Lieb-Luttinger theorems.
The third author wishes to thank William Beckner for his many  constructive comments on the earlier version of the paper posted in the arxiv.org and for his insightful suggestions and encouragement
on possible further study in this direction. His very helpful comments have led to much improvement on its exposition of the paper as well.

\section{The proof of Theorem \ref{theorem1}}
Throughout this section, we shall establish  the Stein-Weiss inequalities with the fractional Poisson kernel. For simplicity, we give the following
notations. Define
$$B_{R}(x)=\{y\in\mathbb{R}^{n}:\ \mid y-x\mid<R,x\in\mathbb{R}^{n}\},$$
$$B_{R}^{n-1}(x)=\{y\in\partial\mathbb{R}^{n}_{+}:\ \mid y-x\mid<R,x\in\partial\mathbb{R}_{+}^{n}\},$$
$$B_{R}^{+}(x)=\{y=(y_{1},y_{2},...,y_{n})\in B_{R}(x):\ y_{n}>0,x\in\mathbb{R}_{+}^{n}\}.$$
For $x=0$, we can also write $B_{R}=B_{R}(0)$,  $B_{R}^{n-1}=B_{R}^{n-1}(0)$,  $B_{R}^{+}=B_{R}^{+}(0)$.
The following lemma will be used in the proof of Theorem \ref{theorem1}. This is an analogue on the half space given in \cite{Dou0} of the result proved in \cite{DHK} on the entire space.

\begin{lemma}\label{lemma2.1}
Let $W(x)$ and $U(\xi)$ be non-negative locally integrable functions defined on  $\mathbb{R}_{+}^{n}$ and $\partial\mathbb{R}_{+}^{n}$  respectively,  for $1<p\leq q<\infty$ and $f\geq0$ on $\partial\mathbb{R}_{+}^{n}$,
\begin{equation}\label{wei}
\left(\int_{\mathbb{R}_{+}^{n}}W(x)\left(\int_{B_{|x|}^{n-1}}f(\xi)d\xi\right)^{q}dx\right)^{\frac{1}{q}}
\leq C_{0}(p,q)\left(\int_{\partial\mathbb{R}_{+}^{n}}f^{p}(\xi)U(\xi)d\xi\right)^{\frac{1}{p}},
\end{equation}
holds if and only if
\begin{equation}\label{inf}
A_{0}=\sup_{R>0}\left\{\left(\int_{|x|\geq R} W(x)dx\right)^{\frac{1}{q}}\left(\int_{|\xi|\leq R}U^{1-p'}(\xi)d\xi\right)^{\frac{1}{p'}}\right\}<\infty.
\end{equation}
On the other hand,
\begin{equation}\label{wei2}
\left(\int_{\mathbb{R}_{+}^{n}}W(x)\left(\int_{\partial\mathbb R_{+}^{n} \setminus B_{|x|}^{n-1}}f(\xi)d\xi\right)^{q}dx\right)^{\frac{1}{q}}
\leq C_{0}(p,q)\left(\int_{\partial\mathbb{R}_{+}^{n}}f^{p}(\xi)U(\xi)d\xi\right)^{\frac{1}{p}},
\end{equation}
holds if and only if
\begin{equation}\label{inf2}
A_{1}=\sup_{R>0}\left\{\left(\int_{|x|\leq R} W(x)dx\right)^{\frac{1}{q}}\left(\int_{|\xi|\geq R}U^{1-p'}(\xi)d\xi\right)^{\frac{1}{p'}}\right\}<\infty.
\end{equation}
\end{lemma}

We now continue with the proof of the Theorem \ref{theorem1}.
\begin{proof}
Without loss of generality, we may assume that $f$ is nonnegative. Define the double weighted integral operator associated with the fractional Poisson kernel
$$P(f)(x)=\int_{\partial\mathbb{R}^n_{+}}P(x,\xi,\gamma)f(\xi)d\xi.$$
It is easy to verify that inequality \eqref{PHLS} is equivalent to the following inequality
\begin{equation*}\label{equ}
||P(f)|x|^{-\beta}||_{L^{q}(\mathbb{R}_{+}^{n})}
\leq C_{n,\alpha,\beta,p,q'}||f|\xi|^{\alpha}||_{L^{p}(\partial\mathbb{R}_{+}^{n})}.
\end{equation*}
Since $q>1$, it follows that
$$\|P(f)|x|^{-\beta}\|^{q}_{L^q(\mathbb{R}^n_{+})}\lesssim P_{1}+P_{2}+P_{3},$$
where
\begin{equation*}\begin{aligned}
&P_{1}=\int_{\mathbb{R}^{n}_{+}}\left(|x|^{-\beta}\int_{B_{|x|/2}^{n-1}}\frac{x_{n}f(\xi)}{|x-\xi|^{n+2-\gamma}}d\xi\right)^{q}dx,\\
&P_{2}=\int_{\mathbb{R}^{n}_{+}}\left(|x|^{-\beta}\int_{B_{2|x|}^{n-1}\setminus B_{|x|/2}^{n-1}}\frac{x_{n}f(\xi)}{|x-\xi|^{n+2-\gamma}}d\xi\right)^{q}dx,\\
&P_{3}=\int_{\mathbb{R}^{n}_{+}}\left(|x|^{-\beta}\int_{\partial\mathbb R^{n}_{+}\setminus B_{2|x|}^{n-1}}\frac{x_{n}f(\xi)}{|x-\xi|^{n+2-\gamma}}d\xi\right)^{q}dx.\\
\end{aligned}\end{equation*}

Thus we just need to show
$$P_{j}\leq C_{n,\alpha,\beta,p,q'}\|f|\xi|^{\alpha}\|^{q}_{L^{p}(\partial\mathbb R_{n}^{+})},\ \ j=1,2,3.$$

First, let us examine $P_{1}$. Since $|\xi|\leq\frac{|x|}{2}$ in this case, we derive that
\begin{equation}\begin{split}
P_{1}\lesssim \int_{\mathbb{R}^{n}_{+}}|x|^{-\beta q-(n+1-\gamma)q}\left(\int_{B_{|x|/2}^{n-1}}f(\xi)d\xi\right)^{q}dx.
\end{split}\end{equation}
Taking  $W(x)=|x|^{-\beta q-(n+1-\gamma)q}$ and  $U(\xi)=|\xi|^{\alpha p}$
in \eqref{wei}, we conclude that
$$P_{1}\leq C_{n,\alpha,\beta,p,q'}\|f|\xi|^{\alpha}\|^{q}_{L^{p}(\partial\mathbb R_{n}^{+})}$$
if $W(x)$ and  $U(\xi)$ satisfy \eqref{inf}.
Indeed, since $\alpha<\frac{n-1}{p'}$, then for any $R>0$, one has
\begin{equation}\begin{split}\label{1}
\int_{|x|\geq R} W(x)dx&=\int_{|x|\geq R}|x|^{-\beta q-(n+1-\gamma)q}dx\\
&=\int_{\partial B_{1}^{+}}d\xi\int_{R}^{\infty}t^{-\beta q-(n+1-\gamma)q+n-1}dt\\
&=C_{1}(n,\lambda,\beta,q)R^{-\beta q-(n+1-\gamma)q+n},
\end{split}\end{equation}
and
\begin{equation}\begin{split}\label{2}
\int_{|\xi|\leq R}U^{1-p'}(\xi)d\xi
&=\int_{|\xi|\leq R}(|\xi|^{\alpha p})^{1-p'}d\xi\\
&=\int_{S^{n-2}}d\eta\int_{0}^{R}r^{\alpha p(1-p')+n-2}dr\\
&=C_{2}(n,\lambda,\alpha,p)R^{\alpha p(1-p')+n-1}.
\end{split}\end{equation}
Combining\eqref{tiaojian}, \eqref{1} and \eqref{2}, we derive that
\begin{equation*}\begin{split}\label{3}
\left(\int_{|x|\geq R} W(x)dx\right)^{\frac{1}{q}}\left(\int_{|y|\leq R}U^{1-p'}(y)dy\right)^{\frac{1}{p'}}
&<C(n,\alpha,\beta,\lambda,p)R^{-\beta -(n+1-\gamma)+\frac{n}{q}+\frac{\alpha p(1-p')+n-1}{p'}}\\
&=C(n,\alpha,\beta,\lambda,p).\\
\end{split}\end{equation*}

Next we estimate $P_{3}$. Since $|\xi|\geq2x$  in this case, it follows that $|\xi-x|\geq\frac{|\xi|}{2}$. Taking $W(x)=|x|^{(-\beta+1)q}$ and $U(\xi)=|\xi|^{(n+2-\gamma+\alpha)p}$ in \eqref{wei2}, we obtain
$$P_{3}\lesssim\int_{\mathbb{R}^{n}_{+}}\left(|x|^{(-\beta+1)q}\int_{\partial\mathbb R^{n}_{+}\setminus B_{2|x|}^{n-1}}f(\xi)|\xi|^{-(n+2-\gamma)}d\xi\right)^{q}dx\leq C(n,\alpha,\beta,\lambda,p)\|f|\xi|^{\alpha}\|^{q}_{L^{p}(\partial\mathbb R_{n}^{+})}$$
if the condition \eqref{inf2} is satisfied. In fact, since $\beta<\frac{n}{q}+1$, then for any $R>0$, there holds
\begin{equation}\begin{split}\label{z3}
\int_{|x|\geq R} W(x)dx&=\int_{|x|\geq R}|x|^{(-\beta+1)q}dx\\
&=\int_{\partial B_{1}^{+}}d\xi\int_{R}^{\infty}t^{(-\beta+1)q+n-1}dt\\
&=C_{1}(n,\lambda,\beta,q)R^{(-\beta+1)q+n},
\end{split}\end{equation}

\begin{equation}\begin{split}\label{z4}
\int_{|\xi|\leq R}U^{1-p'}(\xi)d\xi
&=\int_{|\xi|\leq R}(|\xi|^{(n+2-\gamma+\alpha)p})^{1-p'}d\xi\\
&=\int_{S^{n-2}}d\eta\int_{0}^{R}r^{(n+2-\gamma+\alpha)p(1-p')+n-2}dr\\
&=C_{2}(n,\lambda,\alpha,p)R^{(n+2-\gamma+\alpha)p(1-p')+n-1}.
\end{split}\end{equation}
Combining \eqref{z3} and \eqref{z4}, we verify that condition \eqref{inf2} holds.
\vskip0.3cm

We are left to estimate $P_{2}$.
Since $\frac{|x|}{2}<|\xi|<2|x|$ and $\alpha+\beta\geq0$, it is easy to check
$$|x-\xi|^{\alpha+\beta}<3^{\alpha+\beta}|\xi|^{\alpha+\beta}\leq3^{\alpha+\beta}2^{\beta}|x|^{\beta}|\xi|^{\alpha}.$$
Thus
\begin{equation*}\begin{aligned}
P_{2}&=\int_{\mathbb{R}^{n}_{+}}\left(|x|^{-\beta}\int_{B_{2|x|}^{n-1}\setminus B_{|x|/2}^{n-1}}\frac{x_{n}f(\xi)}{|x-\xi|^{n+2-\gamma}}d\xi\right)^{q}dx\\
&\leq\int_{\mathbb{R}^{n}_{+}}\left(\int_{B_{2|x|}^{n-1}\setminus B_{|x|/2}^{n-1}}\frac{x_{n}f(\xi)|\xi|^{\alpha}}{|x-\xi|^{\alpha+\beta+n+2-\gamma}}d\xi\right)^{q}dx\\
&\leq\int_{\mathbb{R}^{n}_{+}}\left(\int_{\partial\mathbb R^{n}_{+}}\frac{x_{n}f(\xi)|\xi|^{\alpha}}{|x-\xi|^{n+2-(\gamma-\alpha-\beta)}}d\xi\right)^{q}dx.\\
\end{aligned}\end{equation*}
According to the hypotheses of Theorem \ref{theorem1}, we can calculate $2\leq\gamma-\alpha-\beta<n$.  Thanks to the integral inequality \eqref{PHLS1},  we get
$$P_{2}\leq C_{n,\alpha,\beta,p,q'}\|f|\xi|^{\alpha}\|^{q}_{L^{p}(\partial\mathbb R^{n}_{+})}.$$
Thus we accomplish the proof of Theorem \ref{theorem1}.
\end{proof}
\medskip
\section{The proof of Theorem \ref{theorem2}}
In this section, we shall employ the rearrangement inequality and Corollary \ref{corollary2} to investigate the existence of maximizers for the maximizing problem
\begin{equation}\label{mini*}
C_{n,\alpha,\beta,p,q'}:= \sup \{\|V(f)\|_{L^q(\mathbb{R}^n_{+})} : f\geq 0 , \|f\|_{L^p(\partial \mathbb{R}^{n}_{+})}=1\}.
\end{equation}.

For simplicity, we first give some notations. Given $f$ a measurable function on $\partial \mathbb{R}^n_+$, $0<r, s<+\infty$, define the Lorentz norm with indices $r$ and $s$ as
\[
\|f\|_{L^{r,s}(\partial \mathbb{R}^{n}_{+})}=
\begin{cases}
\Big(\int_{0}^{\infty}\big(t^{\frac{1}{r}}f^{*}(t)\big)^{s}\frac{dt}{t}\Big)^{\frac{1}{s}}, & \text{if $s<\infty$,}\\
\sup_{t>0}t^{\frac{1}{p}}f^{*}(t), & \text{if $s=\infty$,}\\
\end{cases}
\]
where $f^{*}(t)$ denotes the decreasing rearrangement of $f$. Now we start our proof.
\vskip0.2cm

\emph{Proof of the theorem \ref{theorem2}:} Assume $\{f_{j}\}_j$ is a maximizing sequence for problem \eqref{mini*}, namely
$$\|f_j\|_{L^p(\partial\mathbb{R}^n_{+})}=1\ \ {\rm and}\ \  \lim_{j\to+\infty}\|V(f_j)\|_{L^q(\mathbb{R}^n_{+})}=C_{n,\alpha,\beta,p,q'}.$$
Since $\alpha$, $\beta\geq0$, by the Riesz rearrangement inequality (see \cite{BLL,LiebLoss}), we obtain
$$\|f_{j}^{*}\|_{L^{p}(\partial\mathbb{R}^n_{+})}=\|f_{j}\|_{L^{p}(\partial\mathbb{R}^n_{+})}=1,\ \
\|V(f_{j})\|^{q}_{L^{q}(\mathbb{R}^n_{+})}\leq\|V(f^{*}_{j})\|^{q}_{L^{q}(\mathbb{R}^n_{+})}.$$
Hence we may assume $\{f_{j}\}_j$ is a nonnegative radially  decreasing sequence.
Given any $f\in L^{p}(\partial\mathbb{R}^{n}_{+})$ and $\lambda>0$, define $f_j^{\lambda}(\xi)=\lambda^{-\frac{n-1}{p}}f(\frac{\xi}{\lambda}).$
It is easy to verify that $$\|f_j^{\lambda}\|_{L^{p}(\partial\mathbb{R}^n_{+})}=\|f_j\|_{L^{p}(\partial\mathbb{R}^n_{+})},\ \ \
\|V(f_j^{\lambda})\|_{L^{q}(\mathbb{R}^n_{+})}=\|V(f_j)\|_{L^{q}(\mathbb{R}^n_{+})}.$$
So is $\{f_j^{\lambda}\}_j$ a maximizing sequence for problem \eqref{mini*}.
Set
$$e_{1}=(1,0,...,0)\in\mathbb{R}^{n-1},\ \  a_{j}=\sup_{\lambda>0}f_{j}^{\lambda}(e_{1})=\sup_{\lambda>0}\lambda^{-\frac{n-1}{p}}f_{j}(\frac{e_{1}}{\lambda}).$$
Direct computations yield that
$$0\leq f_{j}(\xi)\leq a_{j}|\xi|^{-\frac{n-1}{p}}\ \ {\rm and}\ \ \|f_{j}\|_{L^{p,\infty}(\partial\mathbb{R}^n_{+})}\leq w_{n-2}^{\frac{1}{p}}a_{j}.$$
With the help of Corollary \ref{corollary2}, we have
\begin{equation*}\begin{aligned}
\|V(f_{j})\|_{L^{q}(\mathbb{R}^n_{+})}&\lesssim\|f_{j}\|_{L^{p,q}(\partial \mathbb{R}^{n}_{+})}\\
&\lesssim\|f_{j}\|_{L^{p,\infty}}^{1-\frac{p}{q}}\|f_{j}\|^{\frac{p}{q}}_{L^{p}}\\
&\lesssim a_{j}^{1-\frac{p}{q}},
\end{aligned}\end{equation*}
which implies that  $a_{j}\geq c_0$ for some $c_0>0$. Select $\lambda_{j}>0$ satisfying $f_{j}^{\lambda_{j}}(e_{1})\geq c_0$.
Replacing the sequence $\{f_{j}\}_j$ by the new sequence
$\{f_{j}^{\lambda_{j}}\}_j$, still denoted by $\{f_{j}\}_j$, then we obtain that $f_{j}(e_{1})\geq c_0$ for any $j$.
On the other hand, for any $R>0$, we also have
\begin{equation*}\begin{aligned}
v_{n-1}f_{j}^p(R)R^{n-1} &\leq \omega_{n-2} \int_{0}^{R}f_{j}^p(r)r^{n-2}dr\\
&\leq \omega_{n-2} \int_{0}^{+\infty}f_{j}^p(r)r^{n-2}dr\\
&=\int_{\partial\mathbb{R}^n_{+}}f_{j}^p(\xi)d\xi=1,
\end{aligned}\end{equation*}
which implies that
\begin{equation}\label{ab}
0\leq f_{j}(\xi)\leq v_{n-1}^{-\frac{1}{p}}|\xi|^{-\frac{n-1}{p}}.
\end{equation}
Following the Lieb's argument based on the Helly theorem, after passing to a subsequence we may find a nonnegative, radially decreasing function $f$ such that $f_{j}\rightarrow f$ almost everywhere in $\partial\mathbb{R}^{n}_{+}$. It follows that $f(\xi)\geq c_0$ for $|\xi|\leq1$ and $\|f\|_{L^{p}(\partial\mathbb{R}^n_{+})}\leq1$.  By Briesz-Lieb theorem (\cite{BH}), we have
\begin{equation}\label{bf}\begin{aligned}
\lim_{j\rightarrow +\infty}\|f_{j}-f\|^{p}_{L^{p}(\partial\mathbb{R}^n_{+})}&=\lim_{j\rightarrow +\infty}\|f_{j}\|^{p}_{L^{p}(\partial\mathbb{R}^n_{+})}-\|f\|^{p}_{L^{p}(\partial\mathbb{R}^n_{+})}\\
&=1-\|f\|^{p}_{L^{p}(\partial\mathbb{R}^n_{+})}.
\end{aligned}\end{equation}
From \eqref{ab}, we know
\begin{equation}\label{ldct}
V(f_{j})(x)\lesssim|x|^{\beta}\int_{\partial\mathbb{R}^n_{+}}\frac{x_{n}}{|\xi|^{\alpha}(|\xi-x'|^{2}+x_{n}^{2})^{\frac{n+2-\gamma}{2}}}\frac{1}{|\xi|^{\frac{n-1}{p}}}d\xi.
\end{equation}
Since $0\leq\beta<\frac{n}{q}+1$, then the right hand of \eqref{ldct} is finite.  In view of  the dominated convergence theorem, we derive that $\lim\limits_{j\rightarrow +\infty}V(f_{j})(x)= V(f)(x)$ for $x\in \mathbb{R}^n_{+}$. Applying Briesz-Lieb theorem again, one can also derive
\begin{equation*}\begin{aligned}
\lim_{j\rightarrow +\infty}\|V(f_{j})\|^{q}_{{L^{q}}(\mathbb{R}^n_{+})}&=\|V(f)\|^{q}_{{L^{q}}(\mathbb{R}^n_{+})}+\lim_{j\rightarrow +\infty}\|V(f_{j})-V(f)\|^{q}_{L^{q}(\mathbb{R}^n_{+})}\\
&\leq C_{n,\alpha,\beta,p,q'}^q\|f\|^{q}_{L^{p}(\partial\mathbb{R}^n_{+})}+C_{n,\alpha,\beta,p,q'}^q\lim_{j\rightarrow +\infty}\|f_{i}-f\|^{q}_{L^{p}(\partial\mathbb{R}^n_{+})}.\\
\end{aligned}\end{equation*}
Combing this and the inequality \eqref{bf}, we conclude that
$$1\leq \|f\|^{q}_{L^{p}((\partial\mathbb{R}^n_{+})}+\big(1-\|f\|^{p}_{L^{p}((\partial\mathbb{R}^n_{+})}\big)^{\frac{q}{p}}.$$
Since $q>p$ and $f\neq0$, we must have $\|f\|_{L^{p}(\partial\mathbb{R}^n_{+})}=1$. Hence $f_{j}\rightarrow f$ in $L^{p}(\partial\mathbb{R}^n_{+})$ and $f$ is a maximizer of the maximizing problem \eqref{mini*}. Thus we complete the proof of Theorem \ref{theorem2}.

\section{The proof of Theorem \ref{thmx1}}
In this section, we provide the regularity estimate for positive solutions of the integral system \eqref{system3}. We need the following
regularity lifting lemma (see \cite{CL}).
\medskip

Let $V$ be a topological vector space. Suppose there are two extended norms (i.e., the norm
of an element in $V$ might be infinity) defined on $V$,
$$\|\cdot\|_{X},\   \|\cdot\|_{Y}: V\longrightarrow [0, \infty].$$
Let
$$X:=\{f\in V: \|f\|_{X}<\infty \}\, \,{\rm and}\, \, Y:=\{f\in V: \|f\|_{Y}<\infty \}.$$
The operator $T: X\rightarrow Y$ is said to be contracting if for any $f$, $g\in X$, there exists some constant $\eta \in (0,1)$ such that
\begin{equation}\label{con}
\|T(f)-T(g)\|_Y\leq \eta \|f-g\|_X
\end{equation}
and $T$ is said to be shrinking if for any $f\in X$, there exists some constant $\theta \in (0,1)$ such that
\begin{equation}\label{shr}
\|T(f)\|_Y\leq \theta \|f\|_X.
\end{equation}
\begin{lemma}\label{contraction}
Let $T$ be a contraction map from $X$ into itself and from $Y$ into itself. Assume that
 for any $f\in X$, there exists a function $g\in Z:= X \cap Y$ such that
$f=Tf+g \in X.$ Then $f\in Z.$
\end{lemma}
\begin{remark}\label{remxin1}
It is obvious that for a linear operator $T$,  the conditions \eqref{con} and \eqref{shr} are equivalent.
\end{remark}

Now, we start our proof. Denote
\begin{eqnarray*}
u_{a}(\xi)=\left\{
\begin{array}{cl}
&u(\xi),\ \ \ \ \ |u(\xi)|> a\  {\rm or}\  |\xi|>a,\\
&0,\ \ \ \ \ \ \ \ \ {\rm otherwise}.
\end{array}\right.
\end{eqnarray*}
\begin{eqnarray*}
v_{a}(x)=\left\{
\begin{array}{cl}
&v(x),\ \ \ \ \ |v(x)|> a\  {\rm or}\  |x|>a,\\
&0,\ \ \ \ \ \ \ \ \ {\rm otherwise}.
\end{array}\right.
\end{eqnarray*}
$u_{b}(\xi)=u(\xi)-u_{a}(\xi)$ and $v_{b}(x)=v(x)-v_{a}(x)$.
Define the linear operator $T_1$ as
\begin{eqnarray*}
T_{1}(h)(\xi)=\int_{\mathbb{R}^n_{+}} |\xi|^{-\alpha}P(x,\xi,\gamma)v_a^{q_0-1}(x)h(x)|x|^{-\beta}dx,\ \  \xi\in \partial\mathbb{R}^n_{+}
\end{eqnarray*}
and
\begin{eqnarray*}
T_{2}(h)(x)=\int_{\mathbb{R}^n_{+}} |x|^{-\beta}P(x,\xi,\gamma)u_a^{p_0-1}(\xi)h(\xi)|\xi|^{-\alpha}d\xi,\ \  x \in \mathbb{R}^n_{+}.
\end{eqnarray*}

Since $(u,v)\in L^{p_0+1}(\partial\mathbb{R}^n_{+})\times L^{q_0+1}(\mathbb{R}^n_{+})$ is a pair of positive solutions of the integral system \eqref{system3}, then
\begin{equation*}\begin{split}\label{re1}
u(\xi)&=\int_{\mathbb{R}^n_{+}} |\xi|^{-\alpha}P(x,\xi,\gamma)v^{q_0}(x)|x|^{-\beta}dx,\\
&=\int_{\mathbb{R}^n_{+}} |\xi|^{-\alpha}P(x,\xi,\gamma)(v_a+v_b)^{q_0-1}(x)v(x)|x|^{-\beta}dx\\
&=\int_{\mathbb{R}^n_{+}} |\xi|^{-\alpha}P(x,\xi,\gamma)v_a^{q_0-1}(x)v(x)|x|^{-\beta}dx+\int_{\mathbb{R}^n_{+}} |\xi|^{-\alpha}P(x,\xi,\gamma)v_b^{q_0}(x)|x|^{-\beta}dx\\
&=T_1(v)(\xi)+F(\xi)
\end{split}\end{equation*}
and
\begin{equation*}\begin{split}\label{re1}
v(x)&=\int_{\partial\mathbb{R}^n_{+}}|x|^{-\beta}P(x,\xi,\gamma)u^{p_0}(\xi)|\xi|^{-\alpha} d\xi,\\
&=\int_{\partial\mathbb{R}^n_{+}}|x|^{-\beta}P(x,\xi,\gamma)(u_a+u_b)^{p_0-1}(\xi)u(\xi)|\xi|^{-\alpha} d\xi\\
&=\int_{\partial\mathbb{R}^n_{+}} |\xi|^{-\alpha}P(x,\xi,\gamma)u_a^{p_0-1}(\xi)u(\xi)|x|^{-\beta}d\xi+\int_{\partial\mathbb{R}^n_{+}} |\xi|^{-\alpha}P(x,\xi,\gamma)u_b^{p_0}(\xi)|x|^{-\beta}d\xi\\
&=T_2(u)(x)+G(x),
\end{split}\end{equation*}
where $$F(\xi)=\int_{\mathbb{R}^n_{+}} |\xi|^{-\alpha}P(x,\xi,\gamma)v_b^{q_0}(x)|x|^{-\beta}dx,\ \
G(x)=\int_{\partial\mathbb{R}^n_{+}} |\xi|^{-\alpha}P(x,\xi,\gamma)u_b^{p_0}(\xi)|x|^{-\beta}d\xi.$$
Define the operator $T(h_1,h_2)=\big(T_1(h_2), T_2(h_1)\big)$, equip the product space $L^{r}(\partial\mathbb{R}^n_{+})\times L^{s}(\mathbb{R}^n_{+})$ with the norm $\|(h_1,h_2)\|_{r,s}=\|h_1\|_{L^{r}(\partial\mathbb{R}^n_{+})}+\|h_2\|_{L^{s}(\mathbb{R}^n_{+})}.$ It is easy
to see the product space is complete under these norms.\medskip

Observe that $(u,v)$ solves the equation $(u, v) = T (u, v) + (F , G)$. In order to apply the regularity
lifting lemma by contracting operators (Lemma \ref{contraction}), we fix the indices $r$ and $s$ satisfying
\begin{equation}\label{x3}
\frac{1}{s}-\frac{n-1}{n}\frac{1}{r}=\frac{1}{q_0+1}-\frac{n-1}{n}\frac{1}{p_0+1}.
\end{equation}
Note that the interval conditions in Theorem \ref{thmx1} guarantee the existence of such pairs $(r, s)$. Then to arrive at
the conclusion that $(u , v)\in L^{r}(\partial\mathbb{R}^n_{+})\times L^{s}(\mathbb{R}^n_{+})$, we need to verify the following conditions, for sufficiently large $a$. (Here $T$ is
linear, by Remark \ref{remxin1}, we only need to verify that it is shrinking.)
\vskip0.2cm

(i) $T$ is shrinking from  $L^{p_0+1}(\partial\mathbb{R}^n_{+})\times L^{q_0+1}(\mathbb{R}^n_{+})$ to itself.
\vskip0.1cm

(ii) $T$ is shrinking from  $L^{r}(\partial\mathbb{R}^n_{+})\times L^{s}(\mathbb{R}^n_{+})$ to itself.
\vskip0.1cm

(iii) $(F, G)\in L^{p_0+1}(\partial\mathbb{R}^n_{+})\times L^{q_0+1}(\mathbb{R}^n_{+})\cap L^{r}(\partial\mathbb{R}^n_{+})\times L^{s}(\mathbb{R}^n_{+})$.
\vskip0.2cm

We first show that $T$ is shrinking from  $L^{p_0+1}(\partial\mathbb{R}^n_{+})\times L^{q_0+1}(\mathbb{R}^n_{+})$ to itself. Using the integral inequality \eqref{we1} together with the H\"{o}lder inequality, we obtain for $(h_1, h_2)\in L^{p_0+1}(\partial\mathbb{R}^n_{+})\times L^{q_0+1}(\mathbb{R}^n_{+})$,
\begin{equation*}\begin{split}
\|T_1(h_2)\|_{L^{p_0+1}(\partial\mathbb{R}^n_{+})}&\lesssim\|v_a^{q_0-1}\|_{L^{\frac{q_0+1}{q_0-1}}(\mathbb{R}^n_{+})}\|h_2\|_{L^{q_0+1}(\mathbb{R}^n_{+})}\\
&\lesssim\|v_a\|_{L^{q_0+1}(\mathbb{R}^n_{+})}^{q_0-1}\|h_2\|_{L^{q_0+1}(\mathbb{R}^n_{+})}
\end{split}\end{equation*}
and
\begin{equation*}\begin{split}
\|T_2(h_1)\|_{L^{q_0+1}(\mathbb{R}^n_{+})}&\lesssim\|u_a^{p_0-1}\|_{L^{\frac{p_0+1}{p_0-1}}(\partial\mathbb{R}^n_{+})}\|h_1\|_{L^{p_0+1}(\partial\mathbb{R}^n_{+})}\\
&\lesssim\|u_a\|_{L^{p_0+1}(\partial\mathbb{R}^n_{+})}^{p_0-1}\|h_1\|_{L^{p_0+1}(\partial\mathbb{R}^n_{+})}.
\end{split}\end{equation*}
Choosing sufficient large $a$, in view of the integrability $L^{p_0+1}(\partial\mathbb{R}^n_{+})\times L^{q_0+1}(\mathbb{R}^n_{+})$, we derive that
$$\|T(h_1,h_2)\|_{p_0+1,q_0+1}=\|T_1(h_2)\|_{L^{p_0+1}(\partial\mathbb{R}^n_{+})}+\|T_2(h_1)\|_{L^{q_0+1}(\mathbb{R}^n_{+})}\leq \frac{1}{2}\|(h_1, h_2)\|_{p_0+1,q_0+1}.$$ This shows that $T$ is shrinking from  $L^{p_0+1}(\partial\mathbb{R}^n_{+})\times L^{q_0+1}(\mathbb{R}^n_{+})$ to itself.
\medskip

Next, we show that $T$ is shrinking from  $L^{r}(\partial\mathbb{R}^n_{+})\times L^{s}(\mathbb{R}^n_{+})$ to itself.
We use the same tool as we did in (i), that is, the Stein-Weiss inequality with the fractional Poisson kernel in Theorem \ref{theorem1}. Here, we only prove that $\|T_2(h_1)\|_{L^{s}(\mathbb{R}^n_{+})}\leq \frac{1}{2}\|h_1\|_{L^r(\partial\mathbb{R}^n_{+})}$. In fact,
\begin{equation*}\begin{split}
\|T_2(h_1)\|_{L^{s}(\mathbb{R}^n_{+})}&\leq C\|u_a^{p_0-1}h_1\|_{L^{t}(\partial\mathbb{R}^n_{+})}\\
&\leq C\|u_a\|_{L^{p_0+1}(\partial\mathbb{R}^n_{+})}^{p_0-1}\|h_1\|_{L^{r}(\partial\mathbb{R}^n_{+})},
\end{split}\end{equation*}
in which we choose a sufficiently large $a$ that $C\|u_a\|_{L^{p_0+1}(\partial\mathbb{R}^n_{+})}^{p_0-1}\leq \frac{1}{2}$ since $u \in L^{p_0+1}(\partial\mathbb{R}^n_{+})$. Thus, $\|T_2(h_1)\|_{L^{s}(\mathbb{R}^n_{+})}\leq \frac{1}{2}\|h_1\|_{L^{r}(\partial\mathbb{R}^n_{+})}$ for all
$h_1\in L^{r}(\partial\mathbb{R}^n_{+})$. The indices $r$, $s$, $t$ above satisfy
$$\frac{1}{t}=\frac{p_0-1}{p_0+1}+\frac{1}{r}$$ and
\begin{equation*}\begin{split}
\frac{1}{s}&=\frac{\alpha+\beta+2-\gamma}{n}+\frac{n-1}{n}(\frac{p_0-1}{p_0+1}+\frac{1}{r})\\
&=\frac{n-1}{n}\frac{1}{p_0+1}+\frac{1}{q_0+1}-\frac{n-1}{n}+\frac{n-1}{n}(\frac{p_0-1}{p_0+1}+\frac{1}{r})\\
&=\frac{n-1}{n}\frac{1}{r}+\frac{1}{q_0+1}-\frac{1}{p_0+1}\frac{n-1}{n}.
\end{split}\end{equation*}
It is also easy to check that $\alpha\leq \frac{n-1}{t'}$, $\beta<\frac{n}{s}+1$. Similarly we can estimate $T_1(h_2)$ for any $h_2\in L^{s}(\mathbb{R}^n_{+})$. Thus $\|T(h_1,h_2)\|_{r,s}\leq \frac{1}{2}\|(h_1,h_2)\|_{r,s}$. This shows that $T$ is shrinking from  $L^{r}(\partial\mathbb{R}^n_{+})\times L^{s}(\mathbb{R}^n_{+})$ to itself.
\medskip

Finally, we show that $(F, G)\in L^{p_0+1}(\partial\mathbb{R}^n_{+})\times L^{q_0+1}(\mathbb{R}^n_{+})\cap L^{r}(\partial\mathbb{R}^n_{+})\times L^{s}(\mathbb{R}^n_{+})$. It is evident once one notices that $u_b$ and $v_b$ are uniformly bounded by $a$ in $B_a(0)$. Applying the regularity lifting Lemma \ref{contraction}, we finish the proof of Theorem \ref{thmx1}.

\section{The proof of Theorem \ref{thmx2}}
In this section, we give the asymptotic estimate to the double weighted integral system associated with the fractional Poisson kernel.
\medskip

\emph{The Proof of Theorem \ref{thmx2}:} Assume that $(u,v)\in L^{p_0+1}(\partial\mathbb{R}^n_{+})\times L^{q_0+1}(\mathbb{R}^n_{+})$ is a pair of positive solutions of the integral system \eqref{system3}, we only prove that  if $\frac{1}{q_0}-\frac{n+1+\beta-\gamma}{q_0n}>\frac{\beta-1}{n}$, then $$\lim\limits_{|\xi|\rightarrow 0}u(\xi)|\xi|^{\alpha}=\int_{\mathbb{R}^n_{+}}\frac{v^{q_0}(x)x_n}{|x|^{n+2-\gamma+\beta}}dx.$$
We first show that $\int_{\mathbb{R}^n_{+}}\frac{v^{q_0}(x)}{|x|^{n+1-\gamma+\beta}}dx<+\infty$.
Since $u\in  L^{p_0+1}(\partial\mathbb{R}^n_{+})$, so there exists $\xi_0 \in \partial\mathbb{R}^n_{+}$ satisfying $|\xi_0|<\frac{R}{2}$ such that
$u(\xi_0)<\infty$. Then it follows that for any
$$\int_{\mathbb{R}^n_{+}\setminus B_{R}^{+}}\frac{v^{q_0}(x)x_n}{|x|^{n+2-\gamma+\beta}}dx\lesssim \int_{\mathbb{R}^n_{+}\setminus B_{R}^{+}}P(x,\xi_0,\gamma)v^{q_0}(x)|x|^{-\beta}dx\lesssim|\xi_0|^{\alpha}u(\xi_0)<\infty.$$
Observe that $$\int_{\mathbb{R}^n_{+}}\frac{v^{q_0}(x)x_n}{|x|^{n+2-\gamma+\beta}}dx=\int_{B_{R}^{+}}\frac{v^{q_0}(x)x_n}{|x|^{n+2-\gamma+\beta}}dx+
\int_{\mathbb{R}^n_{+}\setminus B_{R}^{+}}\frac{v^{q_0}(x)x_n}{|x|^{n+2-\gamma+\beta}}dx.$$
By the H\"{o}lder inequality, we derive that \begin{equation*}\begin{split}\label{asy1}
\int_{B_{R}^{+}}\frac{v^{q_0}(x)x_n}{|x|^{n+2-\gamma+\beta}}dx\leq \Big(\int_{B_{R}^{+}}(\frac{1}{|x|^{n+1-\gamma+\beta}})^{t'}dx\Big)^{\frac{1}{t'}}
\Big(\int_{B_{R}^{+}}v^{tq_0}(x)dx\Big)^{\frac{1}{t_0}}.
\end{split}\end{equation*}
In order to obtain that  $\int_{\mathbb{R}^n_{+}}\frac{v^{q_0}(x)}{|x|^{n+1-\gamma+\beta}}dx<+\infty$, it suffices to show that
$(n+1-\gamma+\beta)t'<n$ and $$\frac{1}{q_0t}\in \Big(\frac{\beta-1}{n}, \frac{n+1-\gamma+\beta}{n}\Big)\cap \Big(\frac{1}{q_0+1}-\frac{n-1}{n}\frac{1}{p_0+1}+\frac{\alpha}{n}, \frac{1}{q_0+1}-\frac{n-1}{n}\frac{1}{p_0+1}+\frac{n+2-\gamma+\alpha}{n}\Big).$$
According to $(n+1-\gamma+\beta)t'<n$, we have $\frac{1}{q_0t}<\frac{1}{q_0}-\frac{n+1+\beta-\gamma}{q_0n}$. Obviously,
\begin{equation*}\begin{split}
\frac{1}{q_0}-\frac{n+1+\beta-\gamma}{q_0n}&=\frac{1}{q_0+1}-\frac{1}{q_0}(\frac{n-1}{n}\frac{1}{p_0+1}-\frac{\alpha}{n})\\
&> \frac{1}{q_0+1}-\frac{n-1}{n}\frac{1}{p_0+1}+\frac{\alpha}{n}
\end{split}\end{equation*}
since $q_0>1$, $\frac{1}{p_0+1}>\frac{\alpha}{n-1}$ and $\frac{n-1}{n}\frac{1}{p_0+1}+\frac{1}{q_0+1}=\frac{n+\alpha+\beta+1-\gamma}{n}$.
Thus, we are able to choose $t$ such that $(n+1-\gamma+\beta)t'<n$ and $\|v\|_{L^{q_0t}(\mathbb{R}^n_{+})}<\infty$, which leads to   $\int_{\mathbb{R}^n_{+}}\frac{v^{q_0}(x)}{|x|^{n+1-\gamma+\beta}}dx<+\infty$. Next, we prove that
$$\lim\limits_{|\xi|\rightarrow 0}u(\xi)|\xi|^{\alpha}=\int_{\mathbb{R}^n_{+}}\frac{v^{q_0}(x)x_n}{|x|^{n+2-\gamma+\beta}}dx.$$
Observe that
\begin{equation*}\begin{split}
&|\int_{\mathbb{R}^n_{+}}P(x,\xi,\gamma)v^{q_0}(x)|x|^{-\beta}dx-\int_{\mathbb{R}^n_{+}}\frac{v^{q_0}(x)x_n}{|x|^{n+2-\gamma+\beta}}dx\Big|\\
&\ \ \leq \Big|\int_{B_{\delta}^{+}}\Big(P(x,\xi,\gamma)|x|^{-\beta}-\frac{v^{q_0}(x)x_n}{|x|^{n+2-\gamma+\beta}}\Big)dx\Big|+ \Big|\int_{\mathbb{R}^n_{+}
\setminus B_{\delta}^{+}}\Big(P(x,\xi,\gamma)|x|^{-\beta}-\frac{v^{q_0}(x)x_n}{|x|^{n+2-\gamma+\beta}}\Big)dx\Big|\\
&\ \ =J_1+J_2.
\end{split}\end{equation*}
For $J_1$, since \begin{equation*}\begin{split}
\int_{B_{\delta}^{+}}P(x,\xi,\gamma)|x|^{-\beta}dx&\leq \int_{B_{\delta}^{+}(\xi)}\frac{v^{q_0}(x)}{|x-\xi|^{n+1+\beta-\gamma}}dx
+ \int_{B_{\delta}^{+}}\frac{v^{q_0}(x)}{|x|^{n+1+\beta-\gamma}}dx\\
&\leq 2 \|v^q_0\|_{L^{t}(\mathbb{R}^n_{+})} \|\frac{1}{|x|^{n+1+\beta-\gamma}}\|_{L^{t'}({B_{\delta}^{+}})},
\end{split}\end{equation*}
Then $\lim\limits_{\delta\rightarrow 0}\lim\limits_{|\xi|\rightarrow 0}J_1=0$. For $J_2$, using the Lebesgue dominated convergence theorem, we derive that
 $$\lim\limits_{|\xi|\rightarrow 0}\int_{\mathbb{R}^n_{+}
\setminus B_{\delta}^{+}}\Big(P(x,\xi,\gamma)|x|^{-\beta}-\frac{v^{q_0}(x)x_n}{|x|^{n+2-\gamma+\beta}}\Big)dx=0.$$
Combining the above estimate, we derive that
\begin{equation*}\begin{split}
&\lim\limits_{|\xi|\rightarrow 0}\Big|\int_{\mathbb{R}^n_{+}}P(x,\xi,\gamma)v^{q_0}(x)|x|^{-\beta}dx-\int_{\mathbb{R}^n_{+}}\frac{v^{q_0}(x)x_n}{|x|^{n+2-\gamma+\beta}}dx\Big|\\
&=\lim\limits_{\delta\rightarrow 0}\lim\limits_{|\xi|\rightarrow 0}J_1+\lim\limits_{\delta\rightarrow 0}\lim\limits_{|\xi|\rightarrow 0}J_2=0,
\end{split}\end{equation*}
which accomplishes the proof of Theorem \ref{thmx2}.

\section{The proof of Theorem \ref{theorem4}}
In this section, we will establish the necessary condition for positive solutions to the single weighted integral system \eqref{int} with the fractional Poisson kernel.
For $2\leq\gamma<n$, let  $(u,v)\in L^{p_0+1}(|\xi|^{-\alpha}d\xi, \partial\mathbb{R}^n_{+}) \times L^{q_0+1}(|x|^{-\beta}dx, \mathbb{R}^n_{+})$ be a pair of $C^1$ positive solutions of the following integral system
\begin{equation*}\label{integral}\begin{cases}
u(\xi)=\int_{\mathbb{R}^n_{+}} P(x,\xi,\gamma) |x|^{-\beta} v^{q_0}(x)dx,\ \ \xi\in \partial \mathbb{R}^n_{+},\\
v(x)=\int_{\partial\mathbb{R}^n} P(x,\xi,\gamma) |\xi|^{-\alpha} u^{p_0}(\xi)d\xi, \ \ x\in \mathbb{R}^n_{+}.
\end{cases}\end{equation*}

The integration by parts yields that
\begin{equation*}\begin{split}
&\int_{B_R^{n-1}\setminus B_{\varepsilon}^{n-1}}|\xi|^{-\alpha}u^{p_0}(\xi)(\xi\cdot\nabla u(\xi))d\xi\\
&\ \ =\frac{1}{1+p_0}\int_{B_R^{n-1}\setminus B_{\varepsilon}^{n-1}}|\xi|^{-\alpha} \xi\cdot\nabla (u^{1+p_0}(\xi))d\xi\\
&\ \ =\frac{R^{1-\alpha}}{1+p_0}\int_{\partial B_R^{n-1} }u^{p_0+1}(\xi)d\sigma+\frac{\varepsilon^{1-\alpha}}{1+p_0}\int_{\partial B_{\varepsilon}^{n-1} }u^{p_0+1}(\xi)d\sigma\\
&\ \ \ \ -\frac{n-1-\alpha}{1+p_0}\int_{B_R^{n-1}\setminus B_{\varepsilon}^{n-1}} u^{p_0+1}(\xi)|\xi|^{-\alpha}d\xi.\\
\end{split}\end{equation*}

and
\begin{equation*}\begin{split}
&\int_{B_R^{+}\setminus B_{\varepsilon}^{+}}|x|^{-\beta} v^{q_0}(x)(x\cdot\nabla v(x))dx\\
&\ \ =\frac{R^{1-\beta}}{q_0+1}\int_{\partial B_R^{+}}v^{q_0+1}(x)d\sigma+\frac{\varepsilon^{1-\beta}}{q_0+1}\int_{\partial B_\varepsilon^{+}}v^{q_0+1}(x)d\sigma
-\frac{n-\beta}{q_0+1}\int_{B_R^{+}\setminus B_{\varepsilon}^{+}} v^{q_0+1}(x)|x|^{-\beta}dx.\\
\end{split}\end{equation*}
Since  $(u,v)\in L^{p_0+1}(|\xi|^{-\alpha}d\xi, \partial\mathbb{R}^n_{+}) \times L^{q_0+1}(|x|^{-\beta}dx, \mathbb{R}^n_{+})$,
therefore, there exists $R_j\rightarrow+\infty$ and $\varepsilon_j\rightarrow 0$ such that
$$R_j^{1-\alpha}\int_{\partial B_{R_j}^{n-1}} u^{p_0+1}(\xi)d\sigma\rightarrow 0,\ \ \ R_j^{1-\beta}\int_{\partial B_{R_j}^{+}} v^{q_0+1}(x)d\sigma\rightarrow 0$$
and $$\varepsilon_j^{1-\alpha}\int_{\partial B_{\varepsilon_j}^{n-1}} u^{p_0+1}(\xi)d\sigma\rightarrow 0,\ \ \ \varepsilon_j^{1-\beta}\int_{\partial B_{\varepsilon_j}^{+}} v^{q_0+1}(x)d\sigma\rightarrow 0.$$
Combining the above estimate, we have
\begin{equation}\label{poh}\begin{split}
&\int_{\partial\mathbb{R}^n_{+}} |\xi|^{-\alpha}u^{p_0}(\xi)(\xi\cdot\nabla u(\xi))d\xi+\int_{\mathbb{R}^n_{+}}|x|^{-\beta}v^{q_0+1}(x)(x\cdot\nabla v(x))dx\\
&\ \ \ \ =-\frac{n-1-\alpha}{p_0+1}\int_{\partial\mathbb{R}^n_{+}} u^{p_0+1}(\xi)|\xi|^{-\alpha}d\xi-\frac{n-\beta}{q_0+1}\int_{\mathbb{R}^n_{+}}v^{q_0+1}(x)|x|^{-\beta}dx.\\
\end{split}\end{equation}
In view of the single weighted integral system \eqref{int}, direct calculation leads to
\begin{equation}\begin{split}
\nabla u(\xi)\cdot \xi&=\frac{du(\rho x)}{dp}|_{\rho=0}\\
&=-(n+2-\gamma)\int_{\mathbb{R}^n_{+}}P(x,\xi,\gamma)|x-\xi|^{-2}(\xi-x)\cdot \xi |x|^{-\beta}v^{q_0}(x)dx
\end{split}\end{equation}
and
\begin{equation}\begin{split}
\nabla v(x)\cdot x&=\frac{dv(\rho x)}{d\rho}|_{\rho=0}\\
&=-(n+2-\gamma)\int_{\partial\mathbb{R}^n_{+}}P(x,\xi,\gamma)|x-\xi|^{-2}(x-\xi)\cdot x |\xi|^{-\alpha}u^{p_0}(\xi)d\xi\\
&\ \ +\int_{\partial\mathbb{R}^n_{+}}P(x,\xi,\alpha)|\xi|^{-\alpha}u^{p_0}(\xi)d\xi.
\end{split}\end{equation}
Then, it follows that
\begin{equation*}\begin{split}
&\int_{\partial\mathbb{R}^n_{+}} |\xi|^{-\alpha}u^{p_0}(\xi)(\xi\cdot\nabla u(\xi))d\xi+\int_{\mathbb{R}^n_{+}}|x|^{-\beta} v^{q_0}(x)(x\cdot\nabla v(x))dx\\
&\ \ \ \ =-(n+1-\gamma)\int_{\mathbb{R}^n_{+}}\int_{\partial\mathbb{R}^n_{+}}|\xi|^{-\alpha}P(x,\xi,\alpha)u^{p_0}(\xi)v^{q_0}(x)|x|^{-\beta}d\xi dx\\
&\ \ \ \ =-(n+1-\gamma)\int_{\partial\mathbb{R}^n_{+}}|\xi|^{-\alpha}u^{p_0+1}(\xi)d\xi\\
&\ \ \ \ =-(n+1-\gamma)\int_{\mathbb{R}^n_{+}}|x|^{-\beta}v^{q_0+1}(x)dx.
\end{split}\end{equation*}
This together with \eqref{poh} implies that $\frac{n-1-\alpha}{p_0+1}+\frac{n-\beta}{q_0+1}=n+1-\gamma$. Then we accomplish the proof of Theorem \ref{theorem4}.

\end{document}